\documentclass[11pt,twoside]{article}
\usepackage{amsmath,amsbsy,amsfonts}
\usepackage{graphicx,color}
\def\nd{\noindent}
\def\thend{\rule{3mm}{3mm}}
\def\Re{\mathbb{R}}

\newcommand{\ch}[1]{\mbox{\textsl{ch}}(#1)}
\newtheorem{thm}{Theorem}[section]

\newtheorem{lem}{Lemma}[section]
\newtheorem{rmk}{Remark}[section]

\newtheorem{cor}{Corollary}[section]
\newtheorem{definition}{Definition}[section]
\newcommand{\fim}{\hfill\rule{2mm}{2mm}}

\begin{document}

\setlength{\baselineskip}{6.5mm} \setlength{\oddsidemargin}{8mm}
\setlength{\topmargin}{-3mm}
\title{\Large\sf Existence and multiplicity of solutions for a class of quasilinear elliptic field equation on $\mathbb{R}^{N}$}
\author{\sf
 Claudianor O. Alves\thanks{Partially supported by CNPq - Grant 304036/2013-7} \;\; and \;\;  Alan C.B. dos Santos  \\
 Universidade Federal de Campina Grande\\
 Unidade Acad\^emica de Matem\'atica - UAMat\\
 58.429-900 - Campina Grande - PB - Brazil\\
 coalves@mat.ufcg.edu.br \& alan4carlos@gmail.com }

\pretolerance10000
\date{}
\numberwithin{equation}{section} \maketitle

\begin{abstract}
In this paper, we establish existence and multiplicity of solutions for the following class of quasilinear field equation  
\begin{equation*}
-\Delta u+V(x)u-\Delta_{p}u+W'(u)=0, \,\,\, \mbox{in} \,\,\, \mathbb{R}^{N},  \eqno{(P)}
\end{equation*}
where $u=(u_1,u_2,...,u_{N+1})$, $p>N \geq 2$, $W$ is a singular function and $V$ is a positive continuous function.   

\vspace{0.2cm} \noindent \emph{2000 Mathematics Subject Classification} : 35J60, 35A15


\end{abstract}

\section{Introduction}
In this paper, we are interested in the existence of weak solutions for the following class of quasilinear field equation  
\begin{equation*}
-\Delta u+V(x)u-\Delta_{p}u+W'(u)=0, \,\,\, \mbox{in} \,\,\, \mathbb{R}^{N},  \eqno{(P)}
\end{equation*}
where $u=(u_1,u_2,...,u_{N+1})$, $p>N \geq 2$, $V: \mathbb{R}^{N} \to \mathbb{R}$ is a positive continuous function with
$$
V_0=\displaystyle \inf_{x\in\Re^{N}} V(x)>0  \eqno{(V_0)}
$$
and $W: \mathbb{R}^{N+1} \setminus \{\overline{\xi}\} \to \mathbb{R}$ is a $C^1-$function verifying some conditions, which will be fixed later on, and $\overline{\xi} \in \mathbb{R}^{N}$ is a singular point of $W$, that is,
$$
\lim_{x \to \overline{\xi}}|W(x)|=+\infty.
$$  
Moreover, $\Delta{u} = (\Delta{u}_1, ...,\Delta{u}_{N+1})$ and $ \Delta_p{u}$ denotes the $(N+1)$-vector whose $j$-th component is given by $div (|\nabla u|^{p-2}u_j)$.

The motivation of the present paper comes from the seminal papers by Badiale, Benci and D'Aprile \cite{BBDAprile1,BBDAprile2}, whose the existence, multiplicity and concentration of bound states solutions, with one-peak or multi-peak,  have been established for the following class of quasilinear field equation 
\begin{equation} \label{Eq1}
-h^{2}\Delta u+V(x)u-h^{p}\Delta_{p}u+W'(u)=0, \,\,\, \mbox{in} \,\,\, \mathbb{R}^{N}, 
\end{equation}
where $h$ is a positive parameter and $V$ and $W$ are functions verifying some technical conditions, such as: \\
\noindent {\bf Conditions on $W$}: 
  \begin{description}
	\item[$(W_1)$] $W\in C^{1}(\Omega,\Re)$, where $\Omega=\Re^{N+1}\setminus \{\overline{\xi}\}$ for some $\overline{\xi}\in\Re^{N+1}$ with $|\overline{\xi}|=1$;
  \item[$(W_2)$] $W$ is two times differentiable in $0$;
  \item[$(W_3)$] $W(\xi)\geq W(0)=0$ for all $\xi\in\Omega$;
  \item[$(W_4)$] There are $c,\rho>0$ such that 
  $$
  |\xi|<\rho \Longrightarrow W(\xi+\overline{\xi})>c|\xi|^{-q},
  $$
where 
  $$
  \frac{1}{q}=\frac{1}{N}-\frac{1}{p}, \quad N\geq 2, \quad p>N.
  $$
\end{description}

An example of a function $W$ satisfying the above assumptions is the following 
$$
W(\xi)=\frac{|\xi|^{2}}{|\xi-\overline{\xi}|^{q}}, \quad \forall \xi \in \Omega.
$$

\noindent {\bf Conditions on $V$}: Related to $V$, the authors have assumed that
$$
V \in C^{1}(\mathbb{R}^N,\mathbb{R}) \quad \mbox{and} \quad V_0=\inf_{x \in \mathbb{R}^{N}}V(x)>0.
$$
Moreover the following classes of potentials have been considered: 

\noindent {\bf Class 1-} The potential $V$ is coercive, that is, 
$$ 
V(x) \to +\infty \quad \mbox{as} \quad |x| \to +\infty
$$
\noindent {\bf Class 2-} The potential $V$ verifies 
$$
\liminf_{|x| \to +\infty}V(x)>V_0=\inf_{x \in \mathbb{R}^{N}}V(x).
$$
This class of potentials was introduced by Rabinowitz \cite{R12} to study existence of solution for a P.D.E. of the type 
$$
-h^{2}\Delta u +V(x)u=f(u), \quad \mbox{in} \quad \mathbb{R}^{N},
$$
where $h$ is a positive parameter.

\noindent {\bf Class 3-} The potential $V$ has two local isolated minima, that is, there are $r_0,r_1>0$ and $x_0,x_1 \in \mathbb{R}^{N}$ satisfying
$$
\forall x \in \overline{B}_{r_0}(x_0)\setminus \{x_0\} \,:\,V(x)>V(x_0)
$$
and
$$
\forall x \in \overline{B}_{r_1}(x_1)\setminus \{x_1\} \,:\,V(x)>V(x_1)
$$
with $\overline{B}_{r_0}(x_0) \cap \overline{B}_{r_1}(x_1)=\emptyset$. For this class of potential, the result is also true by considering a finite number of local isolated minima for $V$.

The problem (\ref{Eq1}) for $N=3, h=1$ and $V=0$ has been studied in  Benci, D'Avenia, Fortunato and Pisani \cite{BAFP} and Benci, Fortunato and Pisani \cite{BFP1, BFP2}. For related problems with (\ref{Eq1}) involving others classes of potentials, we cite  Benci, Micheletti and Visetti \cite{BMV1,BMV2}, Benci, Fortunato, Masiello and Pisani \cite{BFMP}, D'Aprile \cite{Daprile1,Daprile2,Daprile3,Daprile4},  Visetti \cite{Daniela}, Musso \cite{Musso} and their references.

In general, in the introduction of the above mentioned papers, the reader will find a very nice physical motivation to study (\ref{Eq1}). For example, it is mentioned that this type of problem is related with the study of soliton-like solutions. Moreover, it is also observed that (\ref{Eq1}) appears in the study of the standing wave solutions for the nonlinear Schr\"odinger equation where the presence of a small diffusion parameter $h$ becomes natural.  

Motivated by cited references, we intend to study the existence and multiplicity of solution for $(P)$ for three new classes of potential $V$. Here, we will consider the  following classes: \\

\noindent {\bf Class 4-} The potential $V$ is $\mathbb{Z}^{N}$-periodic, that is,   
$$
V(x+z)=V(x) \,\,\,  \forall x\in\Re^{N} \quad \mbox{and} \quad \forall z\in\mathbb{Z}^{N}.
$$
\noindent {\bf Class 5-} The potential $V$ is asymptotically periodic, that is, there is a $\mathbb{Z}^{N}$-periodic function $V_\infty:\mathbb{R}^{N} \to \mathbb{R}$ such that
$$
	0<\inf_{x \in \mathbb{R}^{N}}V(x) \leq V(x) \leq V_\infty(x) \quad \forall x \in \mathbb{R}^{N}, \eqno{(V_1)}
$$
	and
$$
	|V(x)-V_\infty(x)| \to 0 \quad \mbox{as} \quad |x| \to +\infty. \eqno{(V_2)}
$$
\noindent {\bf Class 6-} The potential $V$ induces a {\it compactness condition}, that is, considering the Hilbert space 
$$
H=\overline{C_{0}^{\infty}(\Re^{N},\Re^{N+1})}^{\|\cdot\|_{H}}
$$
endowed with the norm
$$
\|u\|_{H}=\left(\int_{\Re^{N}}(|\nabla u|^{2}+V(x)|u|^{2})dx \right)^{\frac{1}{2}},
$$
we assume that the embedding $ H \hookrightarrow L^{2}(\mathbb{R}^{N}) $ is compact. 
	
Next, we cite some potentials which belong to Class 6: \\

\noindent 1) \, $V$ is coercive, that is,
$$
V(x) \to +\infty \quad \mbox{as} \quad |x| \to +\infty.
$$

\noindent 2) \, $\displaystyle \frac{1}{V} \in L^{1}(\mathbb{R}^{N})$. \\

\noindent 3) For all $M>0$, we have that
$$
mes(\{x \in \mathbb{R}^{N}\,:\, V(x) \leq M \})<+\infty.
$$
Hereafter, $mes(B)$ denotes the Lebesgue's measure of a mensurable set $B \subset \mathbb{R}^N$. 

The proof that the potentials above belong to Class 6 follows by using the same ideas explored in the papers by Bartsch and Wang \cite{BW}, Costa \cite{Costa}, Kondrat'ev and Shubin \cite{KS} and  Omana and Willem \cite{OW} .

 Here, we will use variational methods to prove our main result, by adapting some ideas explored in \cite{BBDAprile1,BBDAprile2} and their references. For the case where $V$ is periodic, we have proved a new version of the Splitting lemma, see Section 3. For the case where $V$ is  asymptotically periodic or it induces a compactness condition, we have used the Ekeland's variational principle to get a minimizing sequences, which are $(PS)$ sequence, because this type of sequences are better to apply our arguments, for more details see Section 5.  

In order to apply variational methods, we consider the Banach space 
$$
X=\overline{C_{0}^{\infty}(\Re^{N},\Re^{N+1})}^{\|\cdot\|_{X}}
$$
endowed with the norm
$$
\|u\|_{X}=\left(\int_{\Re^{N}}(|\nabla u|^{2}+V(x)|u|^{2})dx \right)^{\frac{1}{2}}+\left(\int_{\Re^{N}}|\nabla u|^{p} dx\right)^{\frac{1}{p}}
$$
and the set
$$
\Lambda=\left\{u\in X : u(x)\neq \overline{\xi}, \,\, \forall x\in\Re^{N}\right\},
$$
which is an open set in $X$.

Using well known arguments, it is possible to prove that the energy functional $E:\Lambda \to \mathbb{R}$ associated with $(P)$, given by
$$
E(u)=\int_{\Re^{N}}\left(\frac{1}{2}(|\nabla u|^{2}+V(x)|u|^{2}) +\frac{1}{p}|\nabla u|^{p}+ W(u)\right) dx
$$
is well defined, $E \in C^{1}(\Lambda, \mathbb{R})$ and
$$
E'(u)v=\int_{\mathbb{R}^{N}}|\nabla u|^{p-2}\nabla u \nabla v\,dx+\int_{\mathbb{R}^{N}}V(x)|u|^{p-2}uv\,dx-\int_{\mathbb{R}^{N}} W'(u)v dx, 
$$
for all $u \in \Lambda$ and $v \in X.$

From the above commentaries, we observe that $u \in \Lambda$ is a weak solution for $(P)$ if, and only if, $u$ is a critical point of $E$. \\

Our main result is the following

\begin{thm}  \label{T1}
Assume that $(W_1)-(W_4)$ hold. Then, \\
\noindent $i)$ If \, $V$ belongs to Class 4 or 6, problem $(P)$ has infinite nontrivial solutions. \\
\noindent $ii)$ If $V$ belongs to Class 5, problem $(P)$ has at least a nontrivial solution. 
\end{thm}

This paper is organized as follows: In Section 2, we fix some notations and prove some preliminary results.  In Section 3, we prove the existence and multiplicity of solutions for the periodic case. In Section 4,  we study the existence of solution for the  asymptotically periodic case, while in Section 5 we show the main result for Class 6.

\section{Preliminary  results}

The results this section will be true assuming on $V$ only $(V_0)$. However, related to function $W$, we will assume the conditions $(W_1)-(W_4)$. The first lemma establishes some important properties involving the space $X$, which will be used very often in this paper.

\begin{lem} \label{space} The following statements hold:\\
\noindent i) $X$ is continuously embedded in $H^{1}(\mathbb{R}^{N},\mathbb{R}^{N+1})$, $W^{1,p}(\mathbb{R}^{N},\mathbb{R}^{N+1})$  and $L^{\infty}(\mathbb{R}^{N},\mathbb{R}^{N+1})$. \\
\noindent ii) For each $u \in X$, 
$$
\lim_{|x| \to +\infty}u(x)=0.
$$
\noindent iii) If $(u_n)$ converges weakly in $X$ to some function $u$, then it converges uniformly on every compact set in $\mathbb{R}^{N}$. 
\end{lem}
\noindent {\bf Proof.} See \cite{BBDAprile1,BBDAprile2}. \fim

\subsection{Topological charge }

In this subsection for convenience of the reader, we repeat the definition of the Topological charge found in \cite{BBDAprile1,BBDAprile2,BFP1,BFP2} and recall some of its main properties. In the open set $\Omega=\mathbb{R}^{N+1} \setminus \{\overline{\xi}\}$, we consider the $N-$sphere centered at $\overline{\xi}$
$$
\Sigma=\left\{ \xi \in \mathbb{R}^{N+1}\,:\; |\xi-\overline{\xi}|=1 \right\}.
$$

On $\Sigma$ we take the north and the south pole, denoted by $\xi_N$ and $\xi_S$, with respect to the axis the origin with $\overline{\xi}$, that is, 
$$
\xi_N=2\overline{\xi} \quad \mbox{and} \quad \xi_S=0.
$$
Then, we consider the projection $P:\Omega \to \Sigma$ defined by 
$$
P(\xi)=\overline{\xi}+\frac{\xi- \overline{\xi}}{|\xi - \overline{\xi}|} \quad \forall \xi \in \Omega.
$$
Notice that, by definition, 
$$
P(\xi)=2\overline{\xi} \Rightarrow |\xi| >1.
$$

\begin{definition}  \label{charge} Given $u \in \Lambda$ and $U \subset \mathbb{R}^{N}$ an open set such that $|u(x)| \leq 1$ on $\partial U$, then we define the (topological) charge of $u$ in the set $U$ as the following integer number
$$
ch(u,U)=deg(P \circ u, U \cap K(u),2\overline{\xi}),
$$
where $K(u)$ is the open set 
$$
K(u)=\left\{ x \in \mathbb{R}^{N}\,:\,|u(x)|>1 \right\}.
$$
Moreover, given $u \in \Lambda$, we define the (topological) charge of $u$ as the integer number
$$
ch(u)=deg(P \circ u,B_R(0),2\overline{\xi})
$$
for all $R>0$ such that $K(u) \subset B_R(0)$.  
\end{definition}

As an immediate consequence of the above definition, we have the lemma below
 
\begin{lem} \label{convergence} Let $(u_n) \subset \Lambda$ and $u \in \Lambda$ such  that
$$
u_n \to u \quad \mbox{uniformly in} \quad  \mathbb{R}^{N}.
$$
Then, there is $n_0 \in \mathbb{N}$ such that
$$
ch(u_n)=ch(u), \quad \forall n \geq n_0.
$$
\end{lem}

For each $q \in \mathbb{Z}$, we set
$$
\Lambda^{q}=\left\{ u \in \Lambda\,:\, ch(u)=q \right\}.
$$
By Lemma \ref{convergence}, each $\Lambda^{q}$ is open in $X$ with
$$
\Lambda= \displaystyle \bigcup_{q \in \mathbb{Z}}\Lambda^{q}
$$
and 
$$
\Lambda^{i}  \cap \Lambda^{j}= \emptyset \quad \mbox{if} \quad i \not= j. 
$$
Using the above notations, we define the open set
$$
\Lambda^{*}= \displaystyle \bigcup_{q \not= 0}\Lambda^{q}.
$$

Using the properties of the Topological charge, it is easy to check that
$$
\partial \Lambda^{*}=\partial \Lambda=\{u \in X \,:\, u(z)=\overline{\xi} \quad \mbox{for some} \quad z \in \mathbb{R}^{N}\}.
$$

\subsection{Technical results}

\begin{lem} \label{L1} For each $\alpha >0$, there exists $\Delta^{*}>0$ such that for every $u \in \Lambda$,
$$
\|u\|_\infty \geq \alpha \Rightarrow E(u) \geq \Delta^{*}.
$$
Hence,  $\displaystyle \inf_{u \in \Lambda^{*}}E(u)>0$.
\end{lem} 

\noindent {\bf Proof.} Arguing by contradiction, if the lemma does not hold, there are $\alpha >0$ and  $(u_n) \subset \Lambda$ such that
$$
E(u_n) \to 0  \quad  \mbox{and} \quad \|u_n\|_\infty \geq \alpha \quad  \forall n \in \mathbb{N}.
$$
Combining the definition of $E$  with the continuous imbedding $X \hookrightarrow (L^{\infty}(\mathbb{R}^{N}))^{N+1}$, the limit $E(u_n) \to 0$ yields
$$
\|u_n\|_\infty \to 0,
$$
which is an absurd.  To conclude the proof, it is enough to recall 
$$
\|u\|_\infty \geq 1, \quad \forall u \in \Lambda^*,
$$
because applying the above argument, there is $\Delta_1>0$ such that
$$
E(u) \geq \Delta_1, \quad \forall u \in \Lambda^*,
$$
implying that $\displaystyle \inf_{u \in \Lambda^{*}}E(u)>0$.

\fim

Arguing as in \cite{BFP2}, we have the following lemmas

\begin{lem} \label{L0} For each $u \in \Lambda$ and for every sequence $(u_n) \subset \Lambda$, if $(u_n)$ weakly converges to $u$, then
$$
\liminf_{n \to +\infty}E(u_n) \geq E(u).
$$
\end{lem}

\begin{lem} Let $(u_n) \subset \Lambda$ be a bounded sequence in $X$ and weakly converging to $u \in \partial \Lambda$. Then,
$$
\int_{\mathbb{R}^{N}}W(u_n)dx \to +\infty.
$$
\end{lem}

As a byproduct of the proof of the last lemma, we deduce the result below
\begin{cor}  \label{C1} For $u \in \partial \Lambda$, 
$$
\int_{\mathbb{R}^{N}}W(u)dx = +\infty.
$$
\end{cor}

By using the previous results, we set the functional  $J: \overline{\Lambda^{*}} \to [0,+\infty]$ given by
$$
J(u)=\left\{
\begin{array}{l}
E(u), \quad u \in \Lambda^{*}, \\
+\infty, \quad u \in \partial \Lambda^{*}.
\end{array}
\right.
$$
Here, we have used the fact that $\partial \Lambda^{*}=\partial \Lambda$.

Using the Lemmas \ref{L0}, \ref{L1} and Corollary \ref{C1}, we derive that $J$ is weakly lower semicontinuity, that is, the lemma below occurs.   

\begin{lem} Let $(u_n) \subset \overline{\Lambda^{*}}$ be a sequence and $u \in X$ such that $u_n \rightharpoonup u$ in $X$. Then
$$
\liminf_{n \to +\infty}J(u_n) \geq J(u).
$$
\end{lem} 

\begin{lem} \label{L3} The functional $J$ is bounded from below on $\overline{\Lambda^{*}}$ and 
$$
0<J_\infty=\inf_{u \in \overline{\Lambda^{*}}}J(u)=\inf_{u \in \Lambda^{*}}E(u).
$$
Moreover, there is $(u_n) \subset \Lambda^{*}$ such that
$$
E(u_n) \to J_\infty \quad \mbox{and} \quad E'(u_n) \to 0 \quad \mbox{as} \quad n \to +\infty.
$$
\end{lem}
\noindent {\bf Proof.} Since $J(w) \geq 0$ for all $w \in \overline{\Lambda^{*}}$, the boundedness from below is immediate. Moreover, recalling that
$$
J(u)=+\infty \quad \mbox{for} \quad u \in \partial \Lambda^{*}=\partial \Lambda,
$$
it follows that
$$
J_\infty=\inf_{u \in \overline{\Lambda^{*}}}J(u)=\inf_{u \in \Lambda^{*}}J(u)=\inf_{u \in \Lambda^{*}}E(u)>0.
$$

The previous study permit us to apply the Ekeland's Variational Principle to get a minimizing sequence verifying  
$$
J(u_n) \to J_\infty \quad \mbox{as} \quad n \to +\infty
$$
and
$$
J(w) \geq J(u_n) - \frac{1}{n}\|w-u_n\|, \,\,\, \forall w \in \overline{\Lambda^{*}} \quad \mbox{and} \quad \forall n \in \mathbb{N}.
$$
The last limit gives that $(u_n) \subset \Lambda^{*}$, and so, 
$$
J(u_n)=E(u_n) \quad \forall n \in \mathbb{N}.
$$
As $\Lambda^{*}$ is a open set in $X$, for each $n$ fixed and $v \in X$, there is $t_n>0$ small enough such that   
$$
u_n+tv \in \Lambda^{*}, \quad \forall t \in [0,t_n).
$$
Hence, $J(u_n+tv)=E(u_n+tv)$ and 
$$
E(u_n+tv) \geq E(u_n) - \frac{t}{n}\|v\|, \quad \forall n \in \mathbb{N}.
$$
Using the fact that $ E \in C^{1}(\Lambda, \mathbb{R})$, the last inequality yields 
$$
E'(u_n)v \geq -\frac{1}{n}\|v\|.
$$
Once $v$ is arbitrary, it follows that
$$
\|E'(u_n)\|\leq \frac{1}{n} \quad \forall n \in \mathbb{N},
$$
finishing the proof. 
\fim

Arguing as above, it is possible to prove the following corollary

\begin{cor} \label{cor3} For each $q \in \mathbb{Z}$, the functional $J$ is bounded from below on $\Lambda^{q}$ and 
$$
0<J^q=\inf_{u \in \overline{\Lambda^{q}}}J(u)=\inf_{u \in \Lambda^{q}}E(u).
$$
Moreover, there is $(u_n) \subset \Lambda^{q}$ such that
$$
E(u_n) \to J^q \quad \mbox{and} \quad E'(u_n) \to 0 \quad \mbox{as} \quad n \to +\infty.
$$
\end{cor}

The next lemma is crucial to prove that weak limit of a $(PS)$ sequence for $E$ is a critical point for $E$. However, since it follows by using well known arguments, we omit its proof.

\begin{lem}   
If $(v_n)$ is a $(PS)_{c}$ sequence for $E$, then for some subsequence,  there is $v \in \Lambda$ verifying
$$
v_n \rightharpoonup v \quad \mbox{in} \quad X
$$
and
$$
\nabla v_n(x) \to \nabla v(x) \quad \mbox{and} \quad v_n(x) \to v(x) \quad \mbox{a.e. in} \quad \mathbb{R}^{N}.
$$
\end{lem}
\section{The periodic potential}

In this section, we will show the existence of solution for $(P)$ by supposing that $V$ is periodic. To begin with, we show a version of the Splitting lemma found in \cite{BBDAprile1,BFP2}, for the case where $V$ is periodic.  

\begin{lem}[Splitting lemma] \label{split}
Let $q \in \mathbb{Z}^*$ and $(u_n)\subset \Lambda^{q}$ be a minimizing sequence for $E$ in $\Lambda^{q}$, that is, 
$$
E(u_n) \to \inf_{u \in \Lambda^q}E(u)=J^{q}.
$$
Then, there are $\overline{u}_1 \in \Lambda^{q}, (z_{n})\subset \mathbb{Z}^{N}$ and $\tilde{R}_{1}>0$ such that, for some subsequence,
\begin{equation}\label{eq(1.1)}
v_{n}:=u_n(\cdot+z_{n})\rightharpoonup \overline{u}_1 \quad \mbox{in} \quad X;
\end{equation}
\begin{equation}\label{eq(1.2)}
\|\overline{u}_1\|_{\infty} \geq 1;
\end{equation}
\begin{equation}\label{eq(1.3)}
|u_n(x)|\leq 1,\quad \forall x\in \Re^{N}\setminus B_{R_1}(z_{n}) \quad \mbox{and} \quad \forall n\in\mathbb{N};
\end{equation}
\begin{equation}\label{eq(1.4)}
\ch{u_n}=\ch{\overline{u}_1};
\end{equation}
\begin{equation}\label{eq(1.5)}
\liminf_{n} E(v_{n})=\liminf_{n}E(u_n)=J^q.
\end{equation}
\end{lem}

\noindent {\bf Proof.} The idea explored in the present proof was inspired in the arguments used in \cite{BFP2}, which has treated the case where the potential $V$ is constant. 
 
Fix $\gamma\in(0,1)$ and  $x_{n}^{1}\in\Re^{N}$ be a maximum point for $|u_n|$. Then, $|u_n(x_{n}^{1})|>1$, otherwise we should have 
$$
|u_n(x)|\leq |u_n(x_{n}^{1})|\leq 1 \quad \forall x\in\Re^{N},
$$
implying that
$$
\|u_n\|_{\infty}\leq 1
$$
which is absurd, because $(u_n)\subset \Lambda^{q}$. Setting 
$$
u_{n}^{1}:=u_n(\cdot+x_n^{1}),
$$
we have 
\begin{equation}\label{eq(1.6')}
\|u_{n}^{1}\|_{\infty}=|u_{n}^{1}(0)|>1.
\end{equation}
As $J^q \in (0, +\infty)$, the sequence $(E(u_n))$ is bounded. Thus, there is  $a>0$ such that 
$$
E(u_n)\leq a \,\,\, \forall n\in\mathbb{N}, 
$$
and so,
$$
E_{0}(u_{n}^{1})\leq E_{V(\cdot+x_{n}^{1})}(u_{n}^{1})=E(u_n)\leq a \,\,\, \forall n\in\mathbb{N},
$$
where
$$
E_0(u)=\int_{\Re^{N}}\left(\frac{1}{2}(|\nabla u|^{2}+V_0|u|^{2}) +\frac{1}{p}|\nabla u|^{p}+ W(u)\right) dx.
$$
Recalling that $E_0$ is coercive on $X$, it follows that $(u_n^{1})$ is bounded in $X$. As $X$ is reflexive, there is $\hat{u}_1\in X $ such that, for some subsequence of $(u_n^{1})$, still denoted by itself, 
$$
u_{n}^{1}\rightharpoonup \hat{u}_1 \quad \mbox{in} \quad X.
$$
Then, by (\ref{eq(1.6')}),
$$
\|\hat{u}_1\|_{\infty}\geq 1.
$$
Now, using the fact that $(E_0(u_{n}^{1}))$ is a bounded sequence and that $(u_{n}^{1})\subset \Lambda$, we deduce that $\hat{u}_1\in\Lambda$. 

In what follows, we fix $R_1>0$ verifying 
\begin{equation}\label{eq(1.80)}
|\hat{u}_1(x)|\leq\gamma \quad x\in\Re^{N}\setminus B_{R_1}(0).
\end{equation}
Now, we will consider two cases: 
\begin{description}
  \item[$(A_1)$] for $n$ large enough  
	$$
  |u_n(x)|\leq 1, \quad \forall x\in \Re^{N}\setminus B_{R_1}(x_{n}^{1}).
  $$
  \item[$(B_1)$] Eventually passing to a subsequence
  $$
  \exists y_n\in\Re^{N}\setminus B_{R_1}(x_{n}^{1}) \quad \mbox{such that} \quad |u_n(y_n)|>1.
  $$

\end{description}

\noindent {\bf Analysis of Case $(A_1)$:\,} 
For each $n\in\mathbb{N}$, there is $z_n\in\mathbb{Z}^{N}$ satisfying
$$
|z_n - x_{n}^{1}|\leq \sqrt{N}, \quad \forall n \in \mathbb{N}.
$$
Fixing $w_n:=z_n - x_{n}^{1}$,  there is $w \in \mathbb{R}^{N}$ with $|w|\leq \sqrt{N}$, such that for some subsequence, $w_n \to w$ in $\Re^{N}$. Setting 
$$
v_{n}:=u_n(\cdot+z_n)=u_{n}^{1}(\cdot+z_n - x_{n}^{1})=u_{n}^{1}(\cdot + w_n) \in \Lambda,
$$
it is easy to see that
\begin{equation}\label{eq(1.9)}
\|v_{n}\|_{\infty}=\|u_{n}^{1}\|_{\infty}=|u_{n}^{1}(0)|=|v_{n}(w_n)|>1.
\end{equation}
Moreover, using the fact that $V$ is $\mathbb{Z}^{N}$- periodic, it follows that
$$
E_0(v_{n})\leq E(v_{n})=E(u_n)\leq a.
$$
Then, there exists $\overline{u}_1 \in\Lambda$ such that, for some subsequence, 
$$
v_{n}\rightharpoonup \overline{u}_1 \quad \mbox{in} \quad X.
$$
The boundedness of $(w_n)$ combined with the fact that $(v_n)$ converges uniformly on every compact set in $\mathbb{R}^{N}$ leads to
$$
v_n(w_n) \to \overline{u}_1(w),
$$ 
from where it follows that $\|\overline{u}_1\|_{\infty}\geq 1$. Moreover, as $(E_0(v_{n}))$ is bounded, it follows  that $\overline{u}_1\in \Lambda$. 

Now, considering $\tilde{R}_1=R_1+\sqrt{N}$, we have that 
$$
B_{\tilde{R}_1}(0)\supset B_{R_1}(0).
$$
Hence, if $x\in \Re^{N}\setminus B_{\tilde{R}_1}(0)$,  $x + z_n \in \Re^{N}\setminus B_{{R}_1}(x_n^1)$, and so, 
$$
|v_n(x)|\leq \gamma \quad \forall x \in\Re^{N}\setminus B_{\tilde{R}_1}(0).
$$
Therefore
$$
|u_n(x)|\leq \gamma \quad \forall x \in\Re^{N}\setminus B_{\tilde{R}_1}(z_n), 
$$
from where it follows that
$$
\ch{u_n} = \ch{u_n, B_{\tilde{R}_1}(z_n)} = \ch{v_{n}, \tilde{B}_{R_1}(0)} = \ch{\overline{u}_1}, \tilde{B}_{R_1}(0)  = \ch{{v}},
$$
showing that $v_n, \overline{u}_1 \in \Lambda^*$. Moreover, by periodicity of $V$, 
$$
E(v_n)=E(u_n) \,\,\, \forall n \in \mathbb{N},
$$ 
showing the lemma. 

\noindent {\bf Analysis of Case $(B_1)$:\,} Next, we will show that this case does not hold. To do that, we will suppose that $(B_1)$ holds. Let $x_n^{2}$ be a maximum point of $|u_n|$ in $\Re^{N}\setminus B_{R_1}(x_{n}^{1})$, which must satisfy $|u_n(x_n^{2})|>1$. Define 
$$
u_{n}^{2}:=u_n(\cdot+x_{n}^{2})
$$
and note that
\begin{equation}\label{eq(1.13')}
\|u_n^{2}\|_{\infty}=|u_{n}^{2}(0)|>1
\end{equation}
with
$$
E_0(u_n^2)\leq E_{V(\cdot+x_n^2)}(u_n^2)=E(u_n)\leq \sup_{n \in \mathbb{N}}E(u_n)<+\infty.
$$
Thereby, $(u_n^2)$ is bounded in $X$, and so, there is $\hat{u}_2\in X$ such that, for some subsequence,
$$
u_{n}^{2}\rightharpoonup \hat{u}_2 \quad \mbox{in} \quad X.
$$
Since $(u_{n}^{2})\subset \Lambda$ and $(E_0(u_{n}^{2}))$ is bounded,  $\hat{u}_2\in \Lambda$. Then, by (\ref{eq(1.13')}), 
\begin{equation}\label{eq(1.14)}
\|\hat{u}_2\|_{\infty}\geq 1.
\end{equation}
Moreover, arguing as in \cite{BFP2}, (\ref{eq(1.80)}) yields  
\begin{equation}\label{eq(1.15)}
|x_n^2-x_n^1|\to+\infty \quad \mbox{as} \quad n\to+\infty.
\end{equation}
Next, fix $\xi_n\in\mathbb{Z}^{N}$ such that
$$
|\xi_n - x_n^2|\leq \sqrt{N} \,\,\, \forall n \in \mathbb{N}.
$$
Setting $\tilde{w}_n:=\xi_n-x_{n}^{2}$, for some subsequence, $\tilde{w}_n\to \tilde{w}\in\Re^{N}$ with $|\tilde{w}|\leq \sqrt{N}$. Considering 
$$
v_n^2:=u_n(\cdot+\xi_n)=u_n^2(\cdot+\xi_n -x_n^{2})=u_{n}^{2}(\cdot + \tilde{w}_n) \in \Lambda,
$$
it follows that
$$
\|v_n^2\|_{\infty}=\|u_n\|_{\infty}=|u_n^2(0)|=|v_{n}^{2}(\tilde{w}_n)|>1
$$
and
$$
E_0(v_n^2)\leq E_{V(\cdot+\xi_n)}(v_n^2)=E(v_n)\leq \sup_{n \in\mathbb{N}}E(u_n). 
$$
The above inequality implies that $(v_n^2)$ is bounded and that there is $\overline{u}_2 \in X$ such that for some subsequence, 
$$
v_n^2\rightharpoonup \overline{u}_2 \quad \mbox{in} \quad X \quad \mbox{and} \quad \|\overline{u}_2\|_{\infty}\geq 1.
$$
Using the fact that,
$$
|\xi_n - z_n|\to+\infty \quad \mbox{as} \quad n\to+\infty,
$$
for each $\eta>0$ and $\rho>0$, we have for $n$ large enough that 
$$
B_{\rho}(\xi_{n})\cap B_{\rho}(z_{n})=\emptyset.
$$
In what follows, we will fix $\rho>0$ verifying 
$$
\int_{B_{\rho}(0)^{c}}\left(\frac{1}{2}(|\nabla \overline{u}_i|^{2}+\|V\|_\infty|\overline{u}_i|^{2}) +\frac{1}{p}|\nabla \overline{u}_i|^{p}+ W(\overline{u}_i)\right) dx <\frac{\eta}{2}, \quad \forall i=1,2.
$$
From this, 
\begin{eqnarray*}
  J^q=\liminf_{n}E(u_n) &\geq& \liminf_{n}(E_{B_{\rho}(z_n)}(u_n) + E_{B_{\rho}(\xi_n)}(u_n)) \\
    &\geq& \liminf_{n} E_{B_{\rho}(z_n)}(u_n)+\liminf_{n}E_{B_{\rho}(\xi_n)}(u_n) \\
    &=& \liminf_{n}E_{B_{\rho}(0)}(v_n) + \liminf_{n}E_{B_{\rho}(0)}(v_n^2) \\
    &\geq& E_{B_{\rho}(0)}(\overline{u}_1) + E_{B_{\rho}(0)}(\overline{u}_2) \\
    &>& E(\overline{u}_1) + E(\overline{u}_2) - \eta.
\end{eqnarray*}
Here, we have denoted by $E_A$ the functional given by
$$
E_{A}(u)=\int_{A}\left(\frac{1}{2}(|\nabla u|^{2}+V(x)|u|^{2}) +\frac{1}{p}|\nabla u|^{p}+ W(u)\right) dx.
$$
As $\eta$ is arbitrary in the last inequality, we can ensure that
\begin{equation}\label{eq(1.17)}
J^q \geq E(\overline{u}_1)+E(\overline{u}_2).
\end{equation}
Next, fix $R_2>0$ verifying 
\begin{equation}\label{eq(1.18)}
|\overline{u}_2(x)|\leq \gamma \quad \forall x\in\Re^{N}\setminus B_{R_2}(0)
\end{equation}
Here, we have again two cases:
\begin{description}
  \item[$(A_2)$] for $n$ large 
  $$
  |u_n(x)|\leq 1 \quad \forall x\in \Re^{N}\setminus(B_{n}^{1} \cup B_{n}^{2})
  $$
  onde $B_{n}^{i}=B_{R_i}(z_n^i), i=1,2$
  \item[$(B_2)$] For some subsequence,
  $$
  \exists x \in \Re^{N} \setminus (B_{n}^{1}\cup B_{n}^{2}) \quad \mbox{such that} \quad |u_n(x)|>1.
  $$
\end{description}
If $(A_2)$ holds, 
\begin{eqnarray*}
  \ch{u_n} &=& \ch{u_n, B_{n}^{1}\cup B_{n}^{2}} \\
    &=& \ch{u_n, B_{n}^{1}} + \ch{u_n, B_{n}^{2}}  \\
    &=& \ch{u_n^1, B_{R_1}(0)} +\ch{u_{n}^{2}, B_{R_2}(0)}  \\
    &=& \ch{\overline{u}_1, B_{R_1}(0)} + \ch{\overline{u}_2, B_{R_2}(0)} \\
    &=& \ch{\overline{u}_1} +\ch{\overline{u}_2}.
\end{eqnarray*}
Then, we can suppose that $\ch{\overline{u}_1}\neq 0$, and so, $\overline{u}_1\in\Lambda^{*}$. Since $\|\overline{u}_2\|_{\infty}\geq 1$, by Lemma \ref{L1}, there is $\Delta^*>0$ verifying $E(\overline{u}_2)\geq \Delta^{*}>0$. Therefore, from (\ref{eq(1.17)}), 
$$
J^q \geq E(\overline{u}_1) + E(\overline{u}_2)\geq J^q +\Delta^{*},
$$
which is an absurd. 

If the case $(B_2)$ occurs, we will consider a maximum point of $|u_n|$ in  $\Re^{N} \setminus (B_n^1 \cup B_{n}^{2})$ and we repeat the same argument used in the case $(B_1)$.   This process terminates in a finite number of steps, because after the $\ell$ steps with $\ell \geq 2$, we have that 
$$
\ell \Delta^{*}\leq \sum_{i=1}^{\ell}E(\overline{u}_i)\leq \liminf_{n}E(u_n)\leq \sup_{n \in \mathbb{N}}E(u_n),
$$
leading to 
$$
2\leq \ell \leq \frac{a}{\Delta^{*}}.
$$
Since the process was finished in the case $(A_\ell)$, it follows that 
$$
|u_n(x)|\leq 1 \quad \forall x\in \Re^{N} \setminus \cup_{i=1}^{\ell} B_{n}^{i}.
$$
Therefore,
$$
0\neq \ch{u_n}=\sum_{i=1}^{\ell} \ch{\overline{u}_i}.
$$
Without lost of  generality, we can assume that $\ch{\overline{u}_1}\neq 0$. This way, $\overline{u}_1\in\Lambda^{*}$ and 
$$
J^q =\liminf_{n} E(u_n) \geq \sum_{i=1}^{\ell}E(\overline{u}_i) = E(\overline{u}_1) + \sum_{i=2}^{\ell}E(\overline{u}_i) \geq J^q + (\ell -1)\Delta^{*}>J^q, 
$$
which is an absurd.  \fim

\vspace{0.5 cm}

\begin{rmk} \label{R1} Here, we would like point out that the above arguments can be used to prove a version of Splitting lemma on $\Lambda^{*}$.

\end{rmk}

Now, we are able to prove our multiplicity result for the periodic case.

\vspace{0.5 cm}

\noindent {\bf Proof of Theorem \ref{T1} - i) (Class 4) } \, Let $q \in \mathbb{Z}^*$ and $(u_n)\subset \Lambda_{q}$ be a minimizing sequence for $E$ in $\Lambda^{q}$. By \textbf{Splitting Lemma}, there exist $\overline{u}_1\in\Lambda^{q}$, $(x_{n}^{1})\subset \mathbb{Z}^{N}$ and $R_1>0$ such that 
\begin{equation}\label{eq(1.19)}
u_n^{1}=u_n(\cdot+x_{n}^{1}) \rightharpoonup \overline{u}_1 \quad \mbox{in} \quad X;
\end{equation}
$$
\|\overline{u}_1\|_{\infty}\geq 1;
$$
$$
|u_n(x)|\leq 1 \quad \forall x\in\Re^{N}\setminus B_{R_1}(x_n^1);
$$
$$
ch(u_n)=ch(\overline{u}_1);
$$
\begin{equation}\label{eq(1.20)}
\liminf_{n} E(u_n^1)=\liminf_{n}E(u_n)=J^q=\inf_{u \in \Lambda^q}J(u).
\end{equation}
Thereby, gathering (\ref{eq(1.19)}), (\ref{eq(1.20)}) and the weakly lower semi-continuity of $E$, we get
$$
J^q=\liminf_{n} E(u_{n}^{1})\geq E(\overline{u}_1)\geq J^q,
$$
implying that
$$
E(\overline{u}_1)=J^q.
$$
Hence, $\overline{u}_1 \in \Lambda^q$ is a nontrivial solution for $(P)$.
\fim

\begin{cor} \label{C2} If $V$ is a periodic function and $(W_1)-(W_4)$ holds, then problem $(P)$ has at least a solution $u \in \Lambda^*$.

\end{cor}

\section{The asymptotically periodic potential}

In this section, we will prove the Theorem \ref{T1} - $ii)$. By Lemma \ref{L3}, there is $(u_n) \subset \Lambda^{*}$ such that
$$
E(u_n) \to J^{*}=\inf\{E(u)\,:\, u \in \Lambda^{*}\} \quad \mbox{and} \quad E'(u_n) \to 0.
$$

We will assume that $V$ is not $\mathbb{Z}^{N}$- periodic. Thereby, by $(V_1)$, there is an open  set $ \Theta \subset \mathbb{R}^{N}$such that
$$
V(x) < V_\infty(x), \,\,\, \forall x \in \Theta.
$$

If $w \in \Lambda^{*}$ denotes the solution found in Theorem \ref{T1}-i) and $E^{\infty}$ denotes the energy functional associated with the periodic problem, the last inequality gives 
$$
J^{*} \leq E(w) < E^{\infty}(w)=J_{\infty}=\inf\{E^{\infty}(u)\,:\, u \in \Lambda^{*}\}.
$$

It is well known that $(u_n)$ is bounded, and so, we can assume that there is $u \in X$ verifying
$$
u_n \rightharpoonup u \quad \mbox{in} \quad X \quad \mbox{and} \quad E'(u)=0.
$$
We claim that $u \not=0$. Indeed, otherwise by  $(V_2)$,
$$
|E(u_n)-E^{\infty}(u_n)| \to 0, 
$$
implying that
$$
E^{\infty}(u_n) \to J^{*}.
$$
Since $J^{*} < J_{\infty}$, the above limit yields 
$$
E^{\infty}(u_n) < J_{\infty} \quad \mbox{for} \quad n \quad \mbox{large enough},
$$
which is a contradiction with definition of $J_{\infty}$. Thereby, $u \not=0$, showing that $E$ has a nontrivial critical point, and so, $(P)$ has a nontrivial solution. \fim

\section{Potential versus compactness}

In this section, we will show the multiplicity of solutions when potential $V$ belongs to Class 6. To this end, the lemma below is a key point in our arguments. 

\begin{lem} \label{estimate}
Let $(u_n)\subset \Lambda^{*}$ be a sequence for $E$ with $(E(u_n))$ bounded. Then, there is $\delta >0$ such that  
$$
|u_n(x)-\overline{\xi}|\geq \delta \,\,\, \forall n \in \mathbb{N} \quad \mbox{and} \quad \forall x \in \mathbb{R}^{N}.
$$
\end{lem}

\noindent {\bf Proof.} If the lemma does not hold, there exists $(y_n) \subset \mathbb{R}^{N}$ verifying 
\begin{equation} \label{E1}
|u_n(y_n)-\overline{\xi}| \to 0 \,\,\, \mbox{as} \,\,\, n \to +\infty.
\end{equation}
Since $(E(u_n))$ is bounded, we derive that $(u_n)$ is bounded in $X$, and so, there is $u \in \Lambda$ such that
$$
u_n \rightharpoonup u \quad  \mbox{in} \quad X.
$$
Moreover, we also have that $(u_n)$ converges uniformly on every compact set contained in $\mathbb{R}^{N}$. This way, we can assume that
$$
|y_n| \to +\infty \quad \mbox{and} \quad |u_n(y_n)|  \geq \frac{3}{4} \quad \forall n \in \mathbb{N}.
$$
Otherwise, for some subsequence, there is $y \in \mathbb{R}^N$ such that
$$
y_n \to y \,\,\, \mbox{as} \,\,\, n \to +\infty,
$$
and so,
\begin{equation} \label{E2}
u_n(y_n) \to u(y).
\end{equation}
From (\ref{E1}) and (\ref{E2}), 
$$
u(y)=\overline{\xi},
$$
which is an absurd, because $u \in \Lambda$.

Fix $\gamma\in(0,\frac{1}{2})$ and  $x_{n}^{1}\in\Re^{N}$ be a maximum point for $|u_n|$. Then, $|u_n(x_{n}^{1})|>1$, otherwise we should have 
$$
|u_n(x)|\leq |u_n(x_{n}^{1})|\leq 1 \quad \forall x\in\Re^{N},
$$
implying that
$$
\|u_n\|_{\infty}\leq 1,
$$
which is absurd, because $(u_n)\subset \Lambda^{*}$. Setting 
$$
u_{n}^{1}:=u_n(\cdot+x_n^{1}),
$$
we have 
\begin{equation}\label{eq(1.6)}
\|u_{n}^{1}\|_{\infty}=|u_{n}^{1}(0)|>1.
\end{equation}
Once $(E(u_n))$ is bounded, there is  $a>0$ such that 
$$
E(u_n)\leq a \,\,\, \forall n\in\mathbb{N}.
$$
Hence,
$$
E_{0}(u_{n}^{1})\leq E_{V(\cdot+x_{n}^{1})}(u_{n}^{1})=E(u_n)\leq a, \,\,\, \forall n\in\mathbb{N}.
$$
Denoting by $X_0$ the Banach space 
$$
X_0=\overline{C_{0}^{\infty}(\Re^{N},\Re^{N+1})}^{\|\cdot\|_{X}}
$$
endowed with the norm
$$
\|u\|_{X}=\left(\int_{\Re^{N}}(|\nabla u|^{2}+V_0|u|^{2})dx \right)^{\frac{1}{2}}+\left(\int_{\Re^{N}}|\nabla u|^{p} dx\right)^{\frac{1}{p}},
$$
it follows that $E_0$ is coercive on $X_0$, and so, $(u_n^{1})$ is bounded in $X_0$. Thereby, as $X_0$ is Reflexive, there is $\overline{u}_1\in X_{0}$ such that, for some subsequence of $(u_n^{1})$, still denoted by itself, 
$$
u_{n}^{1}\rightharpoonup \overline{u}_1 \quad \mbox{in} \quad X_0.
$$
Then, by (\ref{eq(1.6)}), 
$$
\|\overline{u}_1\|_{\infty}\geq 1.
$$
Now, using the fact that $(E_0(u_{n}^{1}))$ is a bounded sequence and that $(u_{n}^{1})\subset \Lambda$, we see that $\overline{u}_1\in\Lambda$. 

In what follows, we fix $R_1>0$ verifying 
\begin{equation}\label{eq(1.8)}
|\overline{u}_1(x)|\leq\gamma \quad x\in\Re^{N}\setminus B_{R_1}(0).
\end{equation}

Let $x_n^{2}$ be a maximum point of $|u_n|$ in $\Re^{N}\setminus B_{R_1}(x_{n}^{1})$, which must satisfy $|u_n(x_n^{2})| \geq \frac{3}{4}$ , because $y_{n_1} \in \Re^{N}\setminus B_{R_1}(x_{n}^{1})$ for some $n_1 \in \mathbb{N}$ . Define 
$$
u_{n}^{2}:=u_n(\cdot+x_{n}^{2})
$$
and note that
\begin{equation}\label{eq(1.13)}
\|u_n^{2}\|_{\infty}=|u_{n}^{2}(0)| \geq \frac{3}{4}
\end{equation}
with
$$
E_0(u_n^2)\leq E_{V(\cdot+x_n^2)}(u_n^2)=E(u_n)\leq \sup_{n \in \mathbb{N}}E(u_n)<+\infty.
$$
Thereby, $(u_n^2)$ is bounded in $X_0$, and so, there is $\overline{u}_2\in X_0$ such that, for some subsequence,
$$
u_{n}^{2}\rightharpoonup \overline{u}_2 \quad \mbox{in} \quad X_0.
$$
Since $(u_{n}^{2})\subset \Lambda_0$ and $(E_0(u_{n}^{2}))$ is bounded,  $\overline{u}_2\in \Lambda$. Then, by (\ref{eq(1.13)}), 
\begin{equation}\label{eq(1.14)}
\|\overline{u}_2\|_{\infty}\geq \frac{3}{4}.
\end{equation}
Moreover, as in \cite{BFP2}, we claim that
\begin{equation}\label{eq(1.15)}
|x_n^2-x_n^1|\to+\infty \quad \mbox{as} \quad n\to+\infty.
\end{equation}
Indeed, considering $z_n=x_n^{2}-x_n^{1}$, and supposing by contradiction that $(z_n)$ is bounded, we can assume that
$$
z_n \to z.
$$
As $|z_n|=|x_n^{2}-x_n^{1}|\geq R_1$, we have that $|z| \geq R_1$. Thus,
$$
|u_1(z)|\leq \gamma < \frac{1}{2}. 
$$
On the other hand, we know that 
$$
\frac{3}{4}\leq |u_n(x_n^{2})|=|u_n^{1}(z_n)|,
$$
then
$$
0<\frac{3}{4}-|\overline{u}_1(z)|\leq |u_n^{1}(z_n)|-|\overline{u}_1(z)|.
$$
Once $u_n^{1}(z_n) \to \overline{u}_1(z)$, taking the limit for $n \to +\infty$, we get a contradiction. 

Next, fix $R_2>0$ verifying 
\begin{equation}\label{eq(1.18)}
|\overline{u}_2(x)|\leq \gamma \quad \forall x\in\Re^{N}\setminus B_{R_2}(0).
\end{equation}

Let $x_n^{3}$ be a maximum point of $|u_n|$ in $\Re^{N}\setminus B_{R_2}(x_{n}^{2})$, which must satisfy $|u_n(x_n^{3})| \geq \frac{3}{4}$, because $y_{n_2} \in \Re^{N}\setminus B_{R_2}(x_{n}^{1})$ for some $n_2 \in \mathbb{N}$ . Define 
$$
u_{n}^{3}:=u_n(\cdot+x_{n}^{3})
$$
and note that
\begin{equation}\label{eq(1.13'')}
\|u_n^{3}\|_{\infty}=|u_{n}^{3}(0)| \geq \frac{3}{4}.
\end{equation}

Arguing as above, we must have 
$$
|x_n^{3}-x_n^{i}| \to +\infty \quad \mbox{for} \quad i=1,2.
$$

Repeating the previous arguments, we find sequences $(x_n^{i}) \subset \mathbb{R}^{N}$ for $i \in \mathbb{N}$ with 
$$
|x_n^{i}-x_n^{j}| \to +\infty \quad \mbox{for} \quad i \not= j,
$$
such that the sequences $u_n^{i}(x)=u_n(x+x_n^{i})$ verify 
$$
u_n^{i} \rightharpoonup \overline{u}_i \in X_0, \quad \overline{u}_i \in \Lambda \quad \mbox{and} \quad \|\overline{u}_i\|_{\infty} \geq \frac{3}{4} \quad \mbox{for} \quad \forall i \in \mathbb{N}.
$$
In what follows, we fix  $k \in \mathbb{N}$ satisfying
\begin{equation} \label{k}
k \Delta^{*} > \sup_{n \in \mathbb{N}}E(u_n)+1,
\end{equation}
where $ \Delta^{*}>0$ is obtained applying the Lemma \ref{L1} for the functional $E_0$, that is, 
$$
E_0(\overline{u}_i) \geq \Delta^{*}, \quad \forall i \in \mathbb{N}. 
$$ 
On the other hand, there is $\rho=\rho(k)>0$ such that
$$
\int_{B_{\rho}(0)^{c}}\left(\frac{1}{2}(|\nabla \overline{u}_i|^{2}+V_0|\overline{u}_i|^{2}) +\frac{1}{p}|\nabla \overline{u}_i|^{p}+ W(\overline{u}_i)\right) dx <\frac{1}{k} \quad \forall i \in \{1,2...,k\}.
$$
Thereby, for $n$ large enough 
$$
B_{\rho}(x_{n}^{i})\cap B_{\rho}(x_{n}^{j})=\emptyset \quad \forall j,i \in \{1,2,...,k\} \quad \mbox{and} \quad i \not= j.
$$
From this, 
\begin{eqnarray*}
  \sup_{n \in \mathbb{N}}E(u_n) \geq \liminf_{n}E_0(u_n) &\geq& \liminf_{n}(\sum_{i=1}^{k}E_{0,B_{\rho}(x_n^{i})}(u_n)) \\
    &\geq& \sum_{i=1}^{k}\liminf_{n} E_{0,B_{\rho}(x_n^{i})}(u_n) \\
    &=& \sum_{i=1}^{k}\liminf_{n} E_{0,B_{\rho}(0)}(u_n^{i}) \\
    &\geq& \sum_{i=1}^{k}\liminf_{n} E_{0,B_{\rho}(0)}(\overline{u}_{i})\\
    &>& \sum_{i=1}^{k}E_0(\overline{u}_{i}) - 1 \geq k \Delta^{*}  - 1,
\end{eqnarray*}
leading to
\begin{equation}\label{eq(1.17)}
k \Delta^{*} < \sup_{n \in \mathbb{N}}E(u_n) +1,
\end{equation}
which contradicts (\ref{k}), finishing the proof of lemma. Here, we have denoted by $E_{0,A}$ the functional given by
$$
E_{0,A}(u)=\int_{A}\left(\frac{1}{2}(|\nabla u|^{2}+V_0|u|^{2}) +\frac{1}{p}|\nabla u|^{p}+ W(u)\right) dx.
$$
\fim

\begin{cor} \label{C2} Let $(u_n)\subset \Lambda^{*}$ be a sequence such that $(E(u_n))$ bounded. Then, there exists $M_1>0$ such that
$$
|W'(u_n)u_n|, W(u_n) \leq M_1 |u_n|^{2}, \,\,\, \forall n \in \mathbb{N}.
$$
\end{cor}
\noindent {\bf Proof.}  By Lemma \ref{estimate}, we know that there is $\delta>0$ such that
$$
|u_n(x) - \overline{\xi}| \geq \delta \quad \forall x \in \mathbb{R}^{N} \quad \mbox{and} \quad \forall n \in \mathbb{N}. 
$$
On the other hand, as $(E(u_n))$ is bounded in $\mathbb{R}$, we derive that $(u_n)$ is also bounded in $X$. Then, $(\|u_n\|_{\infty})$ is bounded, and so, there is $M>0$ such that 
$$
|u_n(x)| \leq M \quad \forall n \in \mathbb{N} \quad \mbox{and} \quad \forall x \in \mathbb{R}^{N}.
$$
From the above study, we have that
$$
u_n(x) \in {\cal U}=\{z \in \mathbb{R}^{N}\,:\, \delta \leq |z-\overline{\xi}| \leq M \} \quad \forall n \in \mathbb{N} \quad \mbox{and} \quad \forall x \in \mathbb{R}^{N}.
$$
Using the hypothesis on $W$, we deduce that there is $M_1>0$ such that
$$
|W(z)|,|W'(z)z| \leq M_1|z|^{2} \quad \forall z \in {\cal U},
$$
and so,
$$
|W(u_n(x))|,|W'(u_n(x))u_n(x)| \leq M_1|u_n(x)|^{2} \quad \forall  n \in \mathbb{N} \quad \mbox{and} \quad \forall x \in \mathbb{R}^{N},
$$
showing the corollary.
\fim

\vspace{0.5 cm}

\noindent {\bf Proof of Theorem \ref{T1}- i) (Class 6): } \, For each $q \in \mathbb{Z}^*$, we know that there is $(u_n) \subset \Lambda^{q}$  such that
$$
E(u_n) \to J^{q} \quad \mbox{and} \quad E'(u_n) \to 0.
$$
Moreover, as $(u_n)$ is bounded in $X$, we can assume that
$$
u_n \rightharpoonup u \quad \mbox{in} \quad X \quad \mbox{and} \quad E'(u)=0.
$$ 
Then,
$$
\int_{\mathbb{R}^{N}}|\nabla u|^{2} dx+\int_{\mathbb{R}^{N}}|\nabla u|^{p} dx+\int_{\mathbb{R}^{N}}V(x)|u|^{2} dx= 
- \int_{\mathbb{R}^{N}}\left\langle \nabla{W}(u),u\right\rangle dx.
$$
On the other hand, the equality $E'(u_n)u_n=o_n(1)$ implies that 
$$
\int_{\mathbb{R}^{N}}|\nabla u_n|^{2} dx+\int_{\mathbb{R}^{N}}|\nabla u_n|^{p} dx+\int_{\mathbb{R}^{N}}V(x)|u_n|^{2} dx= 
- \int_{\mathbb{R}^{N}}\left\langle \nabla{W}(u_n),u_n\right\rangle dx + o_n(1).
$$
The compact embedding $X \hookrightarrow L^{2}(\mathbb{R}^{N})$ together with Corollary \ref{C2} yields 
$$
\int_{\mathbb{R}^{N}}\left\langle \nabla{W}(u_n),u_n\right\rangle dx \to  \int_{\mathbb{R}^{N}}\left\langle \nabla{W}(u),u\right\rangle dx,
$$
from where it follows that
\begin{equation} \label{Limit1}
\lim_{n \to +\infty}\left(\int_{\mathbb{R}^{N}}|\nabla u_n|^{2} dx+\int_{\mathbb{R}^{N}}|\nabla u_n|^{p} dx+\int_{\mathbb{R}^{N}}V(x)|u_n|^{2} dx\right)=\mathcal{A},
\end{equation}
where
$$
\mathcal{A}=\int_{\mathbb{R}^{N}}|\nabla u|^{2} dx+\int_{\mathbb{R}^{N}}|\nabla u|^{p} dx+\int_{\mathbb{R}^{N}}V(x)|u|^{2} dx.
$$
Since for some subsequence
$$
\nabla u_n(x) \to \nabla u(x) \quad \mbox{and} \quad u_n(x) \to u(x) \quad \mbox{a.e. in} \quad \mathbb{R}^{N},
$$
the limit (\ref{Limit1}) gives
$$
u_n \to u \quad \mbox{in} \quad X.
$$

Hence, $u \in \Lambda^{q}$ and $E(u)=J^{q}$, implying that $u$ is a nontrivial critical point of $E$, and so, $u$ is a nontrivial solution of $(P)$. \fim

\end{document}